\documentclass[12pt]{amsart}

\usepackage{amssymb}
\usepackage{bbm}
\usepackage{mathrsfs}

\makeatletter
\def\hsmash{\relax 
  \ifmmode\def\next{\mathpalette\mathhsm@sh}\else\let\next\makehsm@sh
  \fi\next}
\def\makehsm@sh#1{\setbox\z@\hbox{#1}\finhsm@sh}
\def\mathhsm@sh#1#2{\setbox\z@\hbox{$\m@th#1{#2}$}\finhsm@sh}
\def\finhsm@sh{\wd\z@\z@ \box\z@}
\makeatother

\newtheorem{fac}{Fact}[section]
\newtheorem{facs}[fac]{Facts}
\newtheorem{lem}[fac]{Lemma}
\newtheorem{prop}[fac]{Proposition}
\newtheorem{theo}[fac]{Theorem}
\newtheorem{coro}[fac]{Corollary}
\newtheorem{crit}[fac]{Criterion}

\theoremstyle{definition}

\newtheorem{defi}[fac]{Definition}
\newtheorem{const}[fac]{Construction}

\theoremstyle{remark}
\newtheorem{rem}[fac]{Remark}

\newtheorem{ex}[fac]{Example}

\newcommand{\pmodulo}[1]{\nobreak\ifinner\mkern8mu\else\mkern18mu\fi
 (\textup{mod}\,\,#1)}

\newcommand{\br}{ }
\newcommand{\brr}{, }

\newcommand{\Aut}{\mathop{\text{\rm Aut}}\nolimits}
\newcommand{\Div}{\mathop{\text{\rm Div}}\nolimits}
\newcommand{\Gal}{\mathop{\text{\rm Gal}}\nolimits}
\newcommand{\Spec}{\mathop{\text{\rm Spec}}\nolimits}
\newcommand{\Bl}{\mathop{\text{\rm Bl}}\nolimits}
\newcommand{\Br}{\mathop{\text{\rm Br}}\nolimits}

\newcommand{\ev}{\mathop{\text{\rm ev}}\nolimits}

\renewcommand{\div}{\mathop{\text{\rm div}}\nolimits}

\newcommand{\bbA}{{\mathbbm A}}
\newcommand{\bbF}{{\mathbbm F}}
\newcommand{\bbP}{{\mathbbm P}}
\newcommand{\bbQ}{{\mathbbm Q}}
\newcommand{\bbR}{{\mathbbm R}}
\newcommand{\bbZ}{{\mathbbm Z}}

\newcommand{\calA}{{\mathscr{A}}}

\newcommand{\calK}{{\mathscr{K}}}
\newcommand{\calS}{{\mathscr{S}}}

\newcommand{\frakn}{\mathfrak{n}}

\newcommand{\Pb}{{\text{\bf P}}}

\newcounter{abc}
\newenvironment{abc}{\begin{list}{\rm \alph{abc}) }%
{\usecounter{abc} \leftmargin=0.0pt \labelsep=0.0pt %
\listparindent=0.0pt \labelwidth=0.0pt \parsep=\smallskipamount %
\itemsep=0.0pt \topsep=0.0pt \partopsep=\smallskipamount}}{\end{list}}

\newcounter{iii}
\newenvironment{iii}{\begin{list}{\rm \roman{iii}) }%
{\usecounter{iii} \leftmargin=0.0pt \labelsep=0.0pt %
\listparindent=0.0pt \labelwidth=0.0pt \parsep=\smallskipamount%
 \itemsep=0.0pt \topsep=0.0pt \partopsep=\smallskipamount}}{\end{list}}

\def\rightend#1#2{{%
 \leavevmode\nobreak\hskip .5em plus 1fil
 \penalty600 \hskip 0pt plus -1filll
 \vadjust{}\nobreak\hskip 0pt plus 1filll%
 #1\parfillskip=#2\relax \par}}

\def\eop{\ifmmode\rule[-22pt]{0pt}{1pt}\ifinner\tag*{$\square$}\else\eqno{\square}\fi\else\rightend{$\square$}{0pt}\fi}

\title[Brauer-Manin obstruction for degree four del Pezzo surfaces]{On the Brauer-Manin obstruction for \\degree four del Pezzo surfaces}

\begin{document}

\author{J\"org Jahnel}

\address{D\'epartement Mathematik\\ Univ.\ \!Siegen\\ \!Walter-Flex-Str.\ \!3\\ D-57068 \!Siegen\\ \!\mbox{Germany}}
\email{jahnel@mathematik.uni-siegen.de}
\urladdr{http://www.uni-math.gwdg.de/jahnel}

\author{Damaris Schindler}

\address{Hausdorff Center for Mathematics\\ Endenicher Allee 62\\ D-53115 Bonn\\
Germany}
\email{damaris.schindler@hausdorff-center.uni-bonn.de}
\urladdr{http://www.math.uni-bonn.de/people/dschindl}


\date{March~23,~2015.}



\begin{abstract}
We show that, for every integer
$1 \leq d \leq 4$
and every finite set
$S$
of places, there exists a degree
$d$
del Pezzo surface
$X$
over
$\bbQ$
such that
$\Br(X)/\Br(\bbQ) \cong \bbZ/2\bbZ$
and the Brauer-Manin obstruction works exactly at the places in
$S$.
For
$d = 4$,
we prove that in all cases, with the exception of
$S = \{\infty\}$,
this surface may be chosen diagonalizably over
$\bbQ$.
\end{abstract}

\maketitle
\thispagestyle{empty}

\section{Introduction}

A del Pezzo surface is a smooth, proper algebraic
surface~$X$
over a field
$K$
with an ample anti-canonical sheaf
$\calK^{-1}$.
Over~an algebraically
closed field, every del Pezzo surface of degree
$d \leq 7$
is isomorphic to
$\Pb^2$,
blown up in
$(9-d)$
points in general position~\cite[The\-o\-rem~24.4.iii)]{Man}.

According to the adjunction formula, a smooth complete intersection of two quadrics in
$\Pb^4$
is del Pezzo. The converse is true, as~well. For every del Pezzo surface of degree four, its anticanonical image is the complete intersection of two quadrics in
$\Pb^4$~\cite[Theorem~8.6.2]{Do}.\smallskip

For~an arbitrary proper variety
$X$
over
$\bbQ$,
the Brauer-Manin obstruction is a phenomenon that can explain failures of weak approximation or even the Hasse~principle. Its~mechanism works as~follows.

Let
$p$
be any prime number. The~Grothendieck-Brauer group is a contravariant functor from the category of~schemes to the category of abelian groups. In~particular, for an arbitrary scheme
$X$
and a
\mbox{$\bbQ_p$-rational}
point
$x \colon \Spec \bbQ_p \to X$,
there is a restriction~homomorphism
$x^*\colon \Br(X) \to \Br(\bbQ_p) \cong \bbQ/\bbZ$.
For~a Brauer class
$\alpha \in \Br(X)$,
we call
$$\ev_{\alpha,p}\colon X(\bbQ_p) \longrightarrow \bbQ/\bbZ \, , \quad x \mapsto x^*(\alpha) \, ,$$
the local evaluation map, associated
to~$\alpha$.
Analogously, for the real place, there is the local evaluation map
$\ev_{\alpha,\infty} \colon X(\bbR) \to \frac12\bbZ/\bbZ$.

Let~us use the notation
$\Omega$
for the set of all places of
$\bbQ$,
i.e.~for the union of all finite primes together
with~$\infty$.
The~local evaluation maps are continuous with respect to the
\mbox{$p$-adic},
respectively real, topologies
on~$X(\bbQ_\nu)$.
Moreover,~it is well-known that
$\ev_{\alpha,\nu}$
is constant for all but finitely many~places.
Thus,~only adelic points
$x = (x_\nu)_{\nu\in \Omega} \in X(\bbA_\bbQ)$
satisfying
\begin{equation}
\label{evalpha}
\sum_{\nu\in \Omega} \ev_{\alpha,\nu} (x_\nu) = 0 \in \bbQ/\bbZ
\end{equation}
may possibly be approximated by
$\bbQ$-rational
points.\smallskip

\looseness-1
We~say that the Brauer class
$\alpha \in \Br(X)$
{\em works\/}
at a
place~$\nu$
if the local evaluation map
$\ev_{\alpha,\nu}\colon X(\bbQ_\nu) \to \bbQ/\bbZ$
is non-constant. This~is in fact a property of the residue class
of~$\alpha$
in
$\Br(X)/\Br(\bbQ)$.

Observe~that if 
$X(\bbA_\bbQ) \neq \emptyset$
and there exists a Brauer class that works at least at a single place then weak approximation is violated
on~$X$.
On the other hand, if 
$\Br(X)/\Br(\bbQ) \cong \bbZ/2\bbZ$
and a generator works at least at one place then there are adelic points fulfilling (\ref{evalpha}). That~is, the Brauer-Manin obstruction cannot explain a violation of the Hasse principle.\smallskip

The~goal of this paper is to investigate which subsets of
$\Omega$
may occur as the set of places, at which a nontrivial Brauer class works, in the situation of a degree four del Pezzo~surface. Our~first main result is as~follows.

\begin{theo}
\label{theo1}
Let\/~$S \!\subset\! \Omega$
be any finite subset. Then~there exists a degree four del Pezzo surface\/
$X$
over\/
$\bbQ$
having a\/
\mbox{$\bbQ$-rational}
point such that\/
$\Br(X)/\Br(\bbQ) \cong \bbZ/2\bbZ$
and the nontrivial Brauer class works exactly at the places
in\/~$S$.
\end{theo}

In~particular, there is the following example.

\begin{ex}
\label{nurunendl}
Let
$X \subset \Pb^4_\bbQ$
be the degree four del Pezzo surface that is given by the~equations
\begin{eqnarray*}
                T_0T_1 &=& T_2^2 + 7T_3^2 \, , \\
  (T_0-4T_1)(T_0-6T_1) &=& T_2^2 + 7T_4^2 \, .
\end{eqnarray*}
Then~$X$
has a
\mbox{$\bbQ$-rational}
point,
$\Br(X)/\Br(\bbQ) \cong \bbZ/2\bbZ$,
and the nontrivial Brauer class works exactly at the infinite~place.
In~particular, the surface
$X(\bbR)$
has two connected components and
\mbox{$\bbQ$-rational}
points on only one of them.
\end{ex}

More~details about this example are given in Remark~\ref{dist}. The particular case that
$S = \{\infty\}$
is perhaps the most interesting one. Indeed,~there is a relation to Mazur's conjecture \cite[Conjecture~1]{Maz} stating that the closure of
$X(\bbQ)$
with respect to the real topology is equal to a union of connected components. Similar examples for other kinds of surfaces are available in the literature, including singular cubic surfaces~\cite[\S3]{SD1}, conic bundles with five singular fibers \cite[\S3]{Maz}, and others.\smallskip

Recall that over an algebraically closed field two quadratic forms are always simultaneously diagonalizable. We say that a degree four del Pezzo surface is diagonalizable
over~$\bbQ$
if the defining quadratic forms are diagonalizable
over~$\bbQ$.

The~surface from Example~\ref{nurunendl} is not diagonalizable
over~$\bbQ$
but only over
$\bbQ(\sqrt{6})$,
as is easily seen using Fact~\ref{charpol}.b.iii). Somewhat~surprisingly, such a behaviour is necessary at this~point. There~is the following~result.

\begin{theo}
\label{theo2}
Let\/~$X$
be a degree four del Pezzo surface over\/
$\bbQ$
having an adelic point and\/
$\alpha\in \Br(X)$
a Brauer class that works exactly at the infinite~place.
Then\/~$X$
is not diagonalizable
over\/~$\bbQ$.
\end{theo}

Our method of proof uses the fact that diagonal degree four del Pezzo surfaces have nontrivial automorphisms. By functoriality, these operate on
$\Br(X)$,
but the induced operation on
$\Br(X)/\Br(\bbQ)$
turns out to be trivial automatically. Therefore,~every
$\alpha \in \Br(X)/\Br(\bbQ)$
induces a homomorphism
$i_\alpha\colon \Aut'(X) \rightarrow \Br(K)$.
Cf.~Construction \ref{witness} for more~details.

Moreover,~we prove that if
$\alpha \in \Br(X)$
works at
$\infty$
then there is an automorphism
$\sigma \in \Aut(X)$
witnessing this, i.e.\ such that
$i_\alpha(\sigma)$
has a nontrivial component at
$\infty$.
From~this, the claim easily~follows.\smallskip

Our~third main result asserts that, for diagonalizable degree four del Pezzo surfaces, the subset
$\{\infty\}$
is the only exception of this~kind.

\begin{theo}
\label{theo3}
Let\/~$S \subset \Omega$
be a finite subset, different from\/
$\{\infty\}$.
Then there exists a diagonalizable degree four del Pezzo surface\/
$X$
over\/
$\bbQ$
having a\/
\mbox{$\bbQ$-rational}
point such that\/
$\Br(X)/\Br(\bbQ) \cong \bbZ/2\bbZ$
and the nontrivial Brauer class works exactly at the places
in\/~$S$.
\end{theo}

Conjecturally, for degree four del Pezzo surfaces, all failures of weak approximation are due to the Brauer-Manin~obstruction. More~precisely, it is conjectured that
$X(\bbQ)$
is dense in 
$$X(\bbA_\bbQ)^{\Br} := \!\!\!\bigcap\limits_{\alpha \in \Br(X)}\!\!\!\! X(\bbA_\bbQ)^\alpha ,$$
for
$\smash{X(\bbA_\bbQ)^\alpha \subseteq X(\bbA_\bbQ)}$
the subset defined by condition~(\ref{evalpha}) and
$X(\bbA_\bbQ)$
endowed with the product topology induced by the
$\nu$-adic
topologies on
$X(\bbQ_\nu)$.

Due~to work of P.~Salberger and A.\,N.~Skorobogatov \cite[Theorem~0.1]{SSk}, this conjecture is proven under the assumption that
$X$
has a
\mbox{$\bbQ$-rational}
point. In~particular, if
$X$~has
a
\mbox{$\bbQ$-rational}
point then the
\mbox{$\bbQ$-rational}
points
on~$X$
are automatically Zariski dense.

Recall~that all the surfaces provided by Theorem~\ref{theo1} have a
\mbox{$\bbQ$-rational}
point. We~may thus blow up
\mbox{$\bbQ$-rational}
points in general position to obtain del Pezzo surfaces of low degree, thereby unconditionally establishing the following~result.

\begin{theo}
\label{theo4}
Let\/~$S \subset \Omega$
be an arbitrary finite~subset and\/
$d \leq 4$
a positive~integer. Then~there exists a del Pezzo surface\/
$X$
of
degree\/~$d$
over\/~$\bbQ$
having a\/
\mbox{$\bbQ$-rational}
point such that\/
$\Br(X)/\Br(\bbQ) \cong \bbZ/2\bbZ$
and the nontrivial Brauer class works exactly at the places
in\/~$S$.
\end{theo}

It~is well-known that every del Pezzo
surface~$X$
of degree at least five has
$\Br(X)/\Br(\bbQ) = 0$.
One~way to see this is to systematically inspect all possible Galois operations on the exceptional curves in a way analogous to \cite[Chapter~III, 8.21--8.23]{Ja} and to apply~\cite[Proposition 31.3]{Man}. Thus,~Theorem~\ref{theo4} cannot have an analogue for del Pezzo surfaces of higher~degree.

At~least for
$d = 5$
and~$7$,
as well as for
$d=6$
under the additional assumption that
$X$~has
an adelic point, there is also a geometric~argument. Indeed,~these surfaces are birationally equivalent to
$\Pb^2_\bbQ$~\cite[Theorem~2.1]{VA}, cf.~\cite[Theorem~29.4]{Man}.\smallskip

\begin{rem}
For~a degree four del Pezzo surface, the group
$\Br(X)/\Br(\bbQ)$
may be isomorphic to either
$0$,
$\bbZ/2\bbZ$,
or~$(\bbZ/2\bbZ)^2$.
In~the cases of degree
$3$,
$2$,
or~$1$,
there are even more options~\cite[Section~31, Table~3]{Man} and~\cite{SD2}. We~do not know whether the analogue of Theorem~\ref{theo4} is true for a prescribed Brauer~group.
\end{rem}

\section{Brauer classes on degree four del Pezzo surfaces}

The~goal of this section is to gather some facts about degree four del Pezzo surfaces that are necessary for the~following. This includes some results on their Brauer groups and finally leads us to a proof of the assertions made in Example~\ref{nurunendl}. In~other words, we show Theorem~\ref{theo1} under the assumption of Theorem~\ref{theo3}. Unless~a specific choice is made, we work in this section over an arbitrary base
field~$K$ 
of
characteristic~$\neq\! 2$.
Let~us denote
by~$\overline{K}$
an algebraic closure
of~$K$.\medskip 

A del Pezzo surface
$X \subset \Pb^4_K$
of degree four is the base locus of a pencil
$(\mu Q^{(1)} + \nu Q^{(2)})_{(\mu:\nu) \in \Pb^1}$
of quadratic forms in five variables with coefficients in the
field~$K$.
The~generic member of the pencil must be of rank five, as otherwise
$X$
would be a~cone. The~condition that
$\det (\mu Q^{(1)} + \nu Q^{(2)}) = 0$
therefore defines a finite subscheme
$\calS_X \subset \Pb^1_K$
of
degree~$5$.

Choosing a different basis of the pencil yields another embedding
of~$\calS_X$
into the projective line. Thus, one may consider the subscheme
$\calS_X \subset \Pb^1_K$
as an invariant of the surface
$X$~itself.
Moreover,~the definition extends to arbitrary intersections of two quadrics in
$\Pb^4$
that are not~cones.

\begin{facs}
\label{charpol}
\begin{abc}
\item
$X$
is nonsingular if and only if the scheme\/
$\calS_X$
is reduced.
\item
Let~$X\subset \bbP^4$
be a smooth intersection of two quadrics. Then the following state\-ments hold.
\begin{iii}
\item
If\/
$\{s_0,\ldots,s_4\} = \calS_X(\overline{K})$
then the quadratic forms\/
$Q_{s_0}, \ldots, Q_{s_4}$
are exactly of
rank\/~$4$.
\item
The cusps of the cones defined by\/
$Q_{s_i} = 0$,
for\/
$i = 0,\ldots,4$,
are in general linear position
in\/~$\Pb^4$,
i.e.~not contained in any~hyperplane.
\item
$X$
is diagonalizable
over\/
$K$
if and only if\/
$\calS_X$
is split
over\/~$K$.
\end{iii}
\end{abc}\medskip

\noindent
{\bf Proof.}
{\em
These~statements are rather well-known. Proofs may be found, for example, in~\cite{Wi}. Concretely, parts~a) and b.i) are implied by \cite[Proposition~3.26]{Wi}. Furthermore, part b.ii) is \cite[Corollaire 3.29]{Wi}, while part b.iii) is \cite[Corollaire~3.30]{Wi}.%
}%
\eop
\end{facs}

Let\/~$X$
be a degree four del Pezzo surface over a
field~$K$
and assume that there is a
$K$-rational
point
$s \in \calS_X(K)$
as well as that the corresponding degenerate
quadric~$Q_s$
has a
\mbox{$K$-rational}
point, different form the~cusp.
Then~there exist four linearly independent linear forms
$l_1, \ldots, l_4$
such that
$$Q_s = l_1l_2 - (l_3^2 - D l_4^2) \, .$$
Furthermore,
$D$
is the discriminant of
$Q_s$,
considered as a quadratic form in four~variables.

Indeed,~$Q_s$
represents zero nontrivially and it is well-known that such a quadratic form splits off a hyperbolic plane~\cite[Chapitre~4, Proposition~3]{Se}. More~geometrically, one may argue as~follows.
Take~$l_1$
to be a linear form that describes the hyperplane tangent to the cone defined
by~$Q_s$
at a nonsingular
\mbox{$K$-rational}
point. The~restriction
of~$Q_s$
to this hyperplane is of rank two. After scaling, it~may be written in the form
$l_3^2 - D l_4^2$.
Finally,~$Q_s - (l_3^2 - D l_4^2)$
is a quadratic form that vanishes on the hyperplane defined
by~$l_1$
and therefore~splits.

The~four linear forms
$l_1,\ldots,l_4$
must be linearly independent as
$Q_s$
is of rank four. Hence,~the quadratic form
$Q_s$
is equivalent to
$T_0T_1 - (T_2^2 - D T_3^2)$,
which has
discriminant~$D$.\smallskip

The~case most interesting for us is when there are two distinct
$K$-rational
points
$s_1, s_2 \in \calS_X(K)$
and the corresponding degenerate
quadrics~$Q_{s_1}$,
$Q_{s_2}$
have the same discriminant. Then
$X$
may be given by a system of equations in the form
\begin{eqnarray}
\label{SDfam1}
l_{11}l_{12} & = & l_{13}^2 - D l_{14}^2 \, , \\
\label{SDfam2}
l_{21}l_{22} & = & l_{23}^2 - D l_{24}^2 \, .
\end{eqnarray}
For~such surfaces, there is a standard way to write down a Brauer class, which goes back, at least to B.~Birch and Sir Peter Swinnerton-Dyer~\cite{BSD}.

\begin{prop}\medskip
\label{Brauerklasse}
Let\/~$X$
be the degree four del Pezzo surface over a
field\/~$K$,
given by the equations~(\ref{SDfam1},\ref{SDfam2}). Assume that\/
$D$
is a non-square
in\/~$K$
and put\/
$\smash{L := K(\sqrt{D})}$.

\begin{abc}
\item
Then the quaternion algebra (see\/ \cite[Section 15.1]{Pi} for the notation)
$$\textstyle \calA := \big( L(X), \tau, \frac{l_{11}}{l_{21}} \big)$$
over the function field\/
$K(X)$
extends to an Azumaya algebra over the whole
of\/~$X$.
Here,~by
$\tau \in \Gal(L(X)/K(X))$,
we denote the nontrivial~element.
\item
In~the case that\/
$K = \bbQ$,
denote by\/
$\alpha \in \Br(X)$
the Brauer class, defined by the extension
of\/~$\calA$.
Let\/
$\nu$
be any (archimedean or non-archimedean) place
of\/~$\bbQ$.
\begin{iii}
\item
Let\/~$x \in X(\bbQ_\nu)$
be a point and assume that, for some\/
$i,j \in \{1,2\}$,
one has\/
$l_{1i}(x) \neq 0$
and\/~$l_{2j}(x) \neq 0$.
Denote~the corresponding quotient\/
$l_{1i}(x)/l_{2j}(x)$
by\/~$q$.
Then\/
$$\ev_{\alpha,\nu}(x) =
\left\{
\begin{array}{cl}
0       & \text{~if\/~} (q,D)_\nu = 1 \, , \\
\frac12 & \text{~if\/~} (q,D)_\nu = -1 \, ,
\end{array}
\right.
$$
for\/
$(q,D)_\nu$
the Hilbert~symbol.
\item
If\/~$\nu$
is split in\/
$L$
then the local evaluation map\/
$\ev_{\alpha,\nu}$
is constantly~zero.
\end{iii}
\end{abc}\smallskip

\noindent
{\bf Proof.}
{\em
a)
First of all,
$\calA$
is, by construction, a cyclic algebra of degree two. In~particular,
$\calA$
is simple~\cite[Section~15.1, Corollary~d]{Pi}.
Moreover,~$\calA$
is obviously a central
\mbox{$K(X)$-algebra}.

To~prove the extendability assertion, it suffices to show that
$\calA$
extends as an Azumaya algebra over each valuation ring that corresponds to a prime divisor
on~$X$.
Indeed,~this is the classical Theorem of Auslander-Goldman for non-singular surfaces \cite[Proposition~7.4]{AG}, cf.~\cite[Chapter~IV, Theorem~2.16]{Mi}.

To~verify this, we observe that the principal divisor
$\div(l_{11}/l_{21}) \in \Div(X)$
is the norm of a divisor
on~$X_L$.
In~fact, it is the norm of the difference of two prime divisors, the conic, given by
$\smash{l_{11} = l_{13} - \sqrt{D}l_{14} = 0}$,
and the conic, given by
$\smash{l_{21} = l_{23} - \sqrt{D}l_{24} = 0}$.
In~particular,
$\calA$
defines the zero element in
$H^2(\langle\sigma\rangle, \Div(X_L))$.
Under~such circumstances, the extendability
of~$\calA$
over the valuation ring corresponding to an arbitrary prime divisor
on~$X$
is worked out in
\cite[Paragraph 42.2]{Man}.\smallskip

\noindent
b.i)
The~quotients
$$\textstyle \frac{l_{11}}{l_{21}} / \frac{l_{11}}{l_{22}} = \frac{l_{23}^2 - Dl_{24}^2}{l_{21}^2} \, , \quad
\frac{l_{12}}{l_{21}} / \frac{l_{12}}{l_{22}} = \frac{l_{23}^2 - Dl_{24}^2}{l_{21}^2} \, ,
\quad {\rm and} \quad
\frac{l_{11}}{l_{21}} / \frac{l_{12}}{l_{21}} = \frac{l_{13}^2 - Dl_{14}^2}{l_{12}^2}$$
are norms of rational functions
from~$L(X)$.
Therefore,~they define the trivial element of
$H^2(\langle\sigma\rangle, K(X_L)^*) \subseteq \Br K(X)$,
and hence
in~$\Br X$.
In~particular, the four expressions
$l_{1i}/l_{2j}$
define the same Brauer~class.

The~general description of the evaluation map, given in \cite[Paragraph~45.2]{Man}, shows that
$\ev_{\alpha,\nu}(x)$
is equal to
$0$
or
$\frac12$
depending on whether
$q$
is in the image of the norm map
$N_{L_\frakn/\bbQ_\nu}\colon L_\frakn^* \to \bbQ_\nu^*$,
or not, for
$\frakn$
a place of
$L$
lying
above~$\nu$.
This~is exactly what is tested by the Hilbert
symbol~$(q,D)_\nu$.\smallskip

\noindent
ii)
If~$\nu$
is split in
$L$
then the norm map
$\smash{N_{K(X_{L_\frakn})/K(X_{\bbQ_\nu})}\colon K(X_{L_\frakn})^* \to K(X_{\bbQ_\nu})^*}$
is surjective. In~particular,
$l_{11}/l_{21} \in K(X_{\bbQ_\nu})^*$
is the norm of a rational function
on~$X_{L_\frakn}$.
Therefore, it defines the zero class in
$H^2(\langle\sigma\rangle, K(X_{L_\frakn})^*) \subseteq \Br K(X_{\bbQ_\nu})$,
and thus
in~$\Br X_{\bbQ_\nu}$.
To~complete the argument, we note that every
\mbox{$\bbQ_\nu$-rational}
point
$x\colon \Spec \bbQ_\nu \to X$
factors via
$X_{\bbQ_\nu}$.
}
\eop
\end{prop}


In~the following, we will make heavy use of the two facts~below. The~first one recalls the explicit description of the situation when the Brauer group of
$X$
is isomorphic to
$(\bbZ/2\bbZ)^2$. 

\begin{fac}
\label{order4}
Let\/~$X$
be a degree four del Pezzo surface
over a
local or global
field\/~$K$.
In~the local field case, assume that\/
$X(K) \neq \emptyset$,
and in the global field case that\/
$X$
has an adelic point.\smallskip

\noindent
Then\/~$\Br(X)/\Br(K) \cong (\bbZ/2\bbZ)^2$
if and only if\/
$\calS_X$
has three distinct points\/
$s_0$,
$s_1$,
$s_2 \in \calS_X(K)$
so that all three discriminants\/
$D_{s_0}, D_{s_1}, D_{s_2}$
are non-squares
in\/~$K$
and coincide up to square~factors.\smallskip

\noindent
In~this case, representatives of the three nontrivial classes may be obtained as follows. Choose a subset\/
$\{s_i,s_j\} \subset \{s_0,s_1,s_2\}$
of size two. Write\/
$X$
accordingly in the form (\ref{SDfam1},\ref{SDfam2}) and take the corresponding Azumaya algebra as described in Proposition~\ref{Brauerklasse}.\medskip

\noindent
{\bf Proof.}
{\em
This is a well-known fact and a proof may be found, for example, in~\cite[Theorem~3.4]{VAV}. Note that the assumption
on~$X$
implies that, for every closed point
$s \in \calS_X$,
the corresponding
rank-$4$
quadric has a regular point over the residue field of
$s$~\cite[Lemma~5.1]{VAV}.
}
\eop
\end{fac}

\begin{fac}
\label{order2}
Let\/~$X$
be a degree four del Pezzo surface
over a
local or global
field\/~$K$.
In~the local field case, suppose that\/
$X(K) \neq \emptyset$,
and in the global field case that\/
$X$
has an adelic point.\smallskip

\noindent
Assume\/~$X$
to be diagonalizable
over\/~$K$.
Let\/~$D_i \in K$
for\/
$0 \leq i \leq 4$
be the five
rank-$4$
discriminants and assume that
$D_0 = D_1 =: D$.

\begin{abc}
\item
Let\/
$D$
be a non-square
in\/~$K$.
Then the Brauer class\/
$\alpha \in \Br(X)$
described in Proposition~\ref{Brauerklasse} is trivial,
i.e.~$\alpha \in \Br(K)$,
if and only if\/
$D_2$,
$D_3$,
and\/~$D_4$
are all~squares
in\/~$K$.
\item
If the conditions in a) hold or all five discriminants\/
$D_i$
are squares
in\/~$K$,
then one has\/
$\Br(X)/\Br(K) \cong 0$.
\end{abc}\medskip

\noindent
{\bf Proof.}
{\em
a)
This~equivalence statement is established in \cite[Proposition~3.3]{VAV}.\smallskip

\noindent
b)
Fact~\ref{order4} above proves that
$\Br(X)/\Br(K)$
is at most of order~two. If~it were of order exactly two then, by \cite[Theorem~3.4]{VAV}, the nontrivial class could be obtained as described in Proposition~\ref{Brauerklasse}. In~particular, only the case that
$D$
is a non-square remains to be considered. However,~as the other three discriminants are squares, this is exactly the situation in which part~a) proves that the Brauer class is~trivial.
}
\eop
\end{fac}

\begin{rem}
Under the assumptions of Fact \ref{order2}, there is an isomorphism
$$\Br(X)/\Br(K) \stackrel{\cong}{\longleftarrow} \ker \!\big(o\colon(\bbZ/2\bbZ)^5 \rightarrow K^*/(K^*)^2 \big) / T \, ,$$
where the homomorphism
$o$
is given by
$o\colon (a_0,\ldots,a_4) \mapsto \big( D_0^{a_0} \!\cdot\ldots\cdot D_4^{a_4} \bmod (K^*)^2 \big)$
and
$T$
is generated by the vector
$(1,\ldots,1)$
and the standard vectors
$e_i$,
for those
$i \in \{0,\ldots,4\}$
for which
$D_i$
is a perfect~square. Note~that
$D_0 \!\cdot\ldots\cdot D_4$
is a perfect square
in~$K$,
cf.~Fact~\ref{square}.

If~$X$
has a
\mbox{$K$-rational}
point~$x$
not lying on any exceptional curve then this may be seen roughly as~follows. The~blow-up
$\Bl_x(X)$
is a cubic surface with a
\mbox{$K$-rational}
line~$E$.
The~planes through
$E$
equip
$\Bl_x(X)$
with a structure of a conic~bundle. One~can show that the five degenerate conics split exactly over
$\smash{K(\sqrt{D_0}), \ldots, K(\sqrt{D_4})}$.
The~result then follows from the same argument as in \mbox{\cite[Proposition~7.1.1]{Sk}}, cf.~\cite[remarks after Theorem~1.1]{BMS}.

It~requires, however, quite a lot more effort to establish not only an abstract isomorphism, but to prove that the Brauer classes obtained are exactly those expected in view of~Proposition~\ref{Brauerklasse}.
\end{rem}

Once one has an explicit description of the Brauer classes, one needs criteria to understand whether or not they evaluate constantly at a given place. For this the following result turns out to be very useful.

\begin{crit}[A.~V\'arilly-Alvarado and B.~Viray]
\label{VV}
Let\/~$X$
be the degree four del Pezzo surface
over\/~$\bbQ$
given by the equations~(\ref{SDfam1},\ref{SDfam2}). Assume that\/
$D$
is a non-square and let\/
$\alpha \in \Br(X)$
be the Brauer class described in Proposition~\ref{Brauerklasse}.\smallskip

\noindent
Then,~for any place\/
$\nu \neq 2,\infty$
such that the reductions
modulo\/~$\nu$
of the quadratic forms in (\ref{SDfam1}) and (\ref{SDfam2}) both have
rank\/~$4$, 
the local evaluation map\/
$\ev_{\alpha,\nu}$
is~constant.\medskip

\noindent
{\bf Proof.}
{\em
This is \cite[Proposition~5.2]{VAV}.
}
\eop
\end{crit}

\noindent
{\bf Proof of Theorem~\ref{theo1} assuming Theorem~\ref{theo3}.}
Theorem~\ref{theo3} solves the problem for every subset
$S \neq \{\infty\}$.
Thus,~in order to establish Theorem~\ref{theo1}, it suffices to verify the assertions made in Example~\ref{nurunendl}.

For~this, one first checks that
$\calS_X$
has exactly three
\mbox{$\bbQ$-rational}
points, corresponding to the quadratic forms independent of the variable
$T_2$,
$T_3$,
and
$T_4$,
respectively, and a point of degree two that splits over the quadratic
field~$\smash{\bbQ(\sqrt{6})}$.
In~particular,
$X$
is nonsingular.

The~discriminants of the three
\mbox{$\bbQ$-rational}
quadratic forms of
rank~$4$
are, up to square factors,
$1$,
$(-7)$,
and~$(-7)$.
Therefore,~Fact~\ref{order4} shows that
$\Br(X)/\Br(\bbQ)$
is at most of
order~$2$.
On~the other hand, by Proposition~\ref{Brauerklasse}, we have a Brauer class
$\alpha \in \Br(X)$
that is given over the function field
$\bbQ(X)$
as the quaternion algebra
$\smash{\big( \bbQ(\sqrt{-7})(X),\tau,\varphi \big)}$
for
$\smash{\varphi := \frac{T_0-4T_1}{T_1}}$.

Next,~we observe that
$X$
has no real points with
$x_0 = x_1 = 0$.
Moreover,~for an arbitrary real point
$x \in X(\bbR)$
such that
$x_1 \neq 0$,
the equations imply
$x_0/x_1 \geq 0$
and
$\smash{(\frac{x_0}{x_1} - 4)(\frac{x_0}{x_1} - 6) \geq 0}$,
hence
$$x_0/x_1 \in [0,4] \quad {\rm or} \quad x_0/x_1 \geq 6 \, .$$
There~exist real points of both kinds, for example
$(1\!:\!1\!:\!1\!:\!0\!:\!\sqrt{2})$
and
$(8\!:\!1\!:\!1\!:\!1\!:\!1)$.
Since~$(-7) < 0$,
we have that
$(q,-7)_\infty$
is the sign
of~$q$.
Thus,
$\ev_{\alpha,\infty}$
distinguishes the two kinds of real points. In~particular,
$\Br(X)/\Br(\bbQ)$
is indeed of order two and the nontrivial element works at the infinite~place.

It~remains to show that it does not work at any other~place. Criterion~\ref{VV} shows constancy of the evaluation map
$\ev_{\alpha,\nu}$
for all finite places
$\nu \neq 2,7$.
Furthermore,~$\ev_{\alpha,2}$
is constant by Proposition~\ref{Brauerklasse}.b.ii), as the prime
$2$
splits
in~$\smash{\bbQ(\sqrt{-7})}$.

Finally,~for the
prime~$7$,
we argue as~follows. Let
$x \in X(\bbQ_7)$
be any
\mbox{$7$-adic}
point
on~$X$.
Normalize~the coordinates
$x_0,\ldots,x_4$
such that each is a
\mbox{$7$-adic}
integer and at least one is a~unit.
If~$7|x_0$
and
$7|x_1$
then the equations imply that all coordinates must be divisible
by~$7$,
a~contradiction. Hence,~at least one of
$x_0$
and~$x_1$
is a~unit.
Modulo~$7$,
we have
$(\overline{x}_0-4\overline{x}_1)(\overline{x}_0-6\overline{x}_1) = \overline{x}_0\overline{x}_1$
(since both expressions are equal
to~$\overline{x}_2^2$),
and this equation has the solutions
$\overline{x}_0 / \overline{x}_1 = 1, 3$
in
$\bbZ/7\bbZ$.
However, the solution
$\overline{x}_0 / \overline{x}_1 = 3$
is contradictory, as then
$\overline{x}_0\overline{x}_1$
would be a non-square. Consequently,~both
$x_0$
and~$x_1$
must be units and
$$\textstyle \frac{x_0-4x_1}{x_1} \equiv -3 \pmodulo 7 \, ,$$
which implies that
$\smash{\frac{x_0-4x_1}{x_1}}$
is a square
in~$\bbQ_7$.
This~shows
$\smash{\big( \frac{x_0-4x_1}{x_1}, -7 \big)_7 = 1}$
and
$\ev_{\alpha,7}(x) = 0$.
\eop

\begin{rem}
\label{dist}
In~Example~\ref{nurunendl}, weak approximation is disturbed in a rather astonishing~way. The~smooth manifold
$X(\bbR)$
is disconnected into two~components. There~are two kinds of real points
$x \in X(\bbR)$,
those with
$x_0/x_1 \in [0,4]$
and those such that
$x_0/x_1 \in [6,\infty]$.
However,~for every
\mbox{$\bbQ$-rational}
point
$x \in X(\bbQ)$,
one has
$x_0/x_1 > 6$.

A naively implemented point search shows that there are exactly 792
\mbox{$\bbQ$-rational}
points of naive height up to
$1000$ on~$X$.
The~smallest value of the quotient
$x_0/x_1$
is~$319/53 \approx 6.019$.
\end{rem}

\section{Diagonal degree four del Pezzo surfaces}

The~goal of this section is to collect some facts about diagonal degree four del Pezzo surfaces. These~will lead us to a proof of Theorem~\ref{theo2}.\medskip

Let~$X$
be a diagonal degree four del Pezzo surface over a base field
$K$,
i.e.~one that is given by equations of the~form
\begin{eqnarray}
\label{diag1}
a_0T_0^2 + \ldots + a_4T_4^2 & = & 0 \, , \\
\label{diag2}
b_0T_0^2 + \ldots + b_4T_4^2 & = & 0
\end{eqnarray}
with coefficients
in~$K$.
Then,~ ~for every
$(i_0,\ldots,i_4) \in \{0,1\}^5$,
the map
$$(T_0:\ldots:T_4) \mapsto ((-1)^{i_0}T_0:\ldots:(-1)^{i_4}T_4)$$
defines a
\mbox{$K$-automorphism}
of~$X$.
Thus,~there is a subgroup
$\Aut'(X) \subseteq \Aut_K(X)$
that is isomorphic to
$(\bbZ/2\bbZ)^4$.

It~is known that the automorphism group of a degree four del Pezzo surface over an algebraically closed field is generically isomorphic to
$(\bbZ/2\bbZ)^4$
and that there are particular cases, where the automorphism group is larger~\cite[Theorem~8.6.8]{Do}.

\begin{lem}
Let\/~$X$
be a diagonal degree four del Pezzo surface over a local or global
field\/~$K$.
In~the local field case, suppose that\/
$X(K) \neq \emptyset$,
and in the global field case that\/
$X$
has an adelic point.\smallskip

\noindent
Then the natural operation of\/
$\Aut'(X)$
on\/
$\Br(X)$
induces the trivial operation on\/
$\Br(X)/\Br(K)$.\medskip

\noindent
{\bf Proof.}
{\em
This is trivially true if
$\Br(X)/\Br(K) \cong 0$
or~$\bbZ/2\bbZ$.
Otherwise,~it follows from the description of the representatives given in Fact~\ref{order4}.
}
\eop
\end{lem}

\begin{const}
\label{witness}
Let~$X$
be a diagonal degree four del Pezzo surface over a local or global
field~$K$.
In~the local field case, suppose that
$X(K) \neq \emptyset$,
and in the global field case that
$X$
has an adelic point.

By functoriality, the operation of
$\Aut'(X)$
on~$X$
induces an operation
on~$\Br(X)$,
which is necessarily trivial
on~$\Br(X)/\Br(K)$.
Thus,~for every
$\alpha \in \Br(X) / \Br(K)$,
there is a natural homomorphism
$$i_\alpha\colon \Aut'(X) \longrightarrow \Br(K) \, ,$$
given by the condition that
$\sigma^* \alpha = \alpha + i(\sigma)$
for
$\sigma \in \Aut'(X)$.
\end{const}

\begin{defi}
Let~$K=\bbQ$
and assume that a Brauer class
$i_\alpha(\sigma)$
in the image of
$i_\alpha$
has a nontrivial component at the
place~$\nu$.
Then,~as
$$\ev_{\alpha,\nu}(\sigma(x)) = \ev_{\alpha,\nu}(x) + i(\sigma)_\nu \, ,$$
the Brauer class certainly works
at~$\nu$.
We~say in this situation that
$\sigma$
is a {\em witness\/} for the non-constancy of the local evaluation map
at~$\nu$.
\end{defi}

\begin{fac}
\label{square}
Let\/~$X$
be a diagonal degree four del Pezzo surface over a
field\/~$K$
and\/
$D_0,\ldots,D_4$
be the discriminants of the five associated quadratic forms of
rank\/~$4$.
Then\/~$D_0\!\cdot\ldots\cdot D_4$
is a square
in\/~$K$.\medskip

\noindent
{\bf Proof.}
{\em
This is a direct calculation.
}
\eop
\end{fac}

\begin{lem}
\label{XR}
Let\/~$X$
be a diagonal degree four del Pezzo surface
over\/~$\bbR$
that has a real~point. Assume~that\/
$\Br(X)/\Br(\bbR) \neq 0$.\smallskip

\noindent
Then\/~$X(\bbR)$
splits into two connected components. Moreover,~there is an element\/
$\sigma \in \Aut'(X)$
that interchanges these~components.\medskip

\noindent
{\bf Proof.}
{\em
By~Fact~\ref{square}, there are three~cases. Either,~no of the three
rank-$4$
discriminants is negative, or exactly two, or exactly four of~them. Fact~\ref{order2}.b) shows that
$\Br(X)/\Br(\bbR) \neq 0$
is possible only in the last~case.

Then~the pencil of quadrics
in~$\Pb^4$
associated with
$X$
contains four 
rank-$4$
quadrics of negative discriminant. We may write each of them in the shape
$$-c_0 T_{i_0}^2 + c_1 T_{i_1}^2 + c_2 T_{i_2}^2 + c_3 T_{i_3}^2 = 0 \, ,$$
for
$c_0, \ldots, c_3 > 0$,
and say that the variable
$T_{i_0}$
is distinguished by the form~considered.

We~claim that not all four forms may distinguish the same variable. Indeed,~if that would be the case then we also had
$-c'_0 T_{i_0}^2 + c'_1 T_{i_1}^2 + c'_2 T_{i_2}^2 + c'_4 T_{i_4}^2 = 0$,
which shows that the form in the pencil that does not involve
$T_{i_0}$
has opposite signs at
$T_{i_3}^2$
and~$T_{i_4}^2$.
The~same argument for all combinations of two of the four quadratic forms enforces six opposite signs among the four coefficients
of~$T_{i_1}^2, \ldots, T_{i_4}^2$,
a~contradiction.

Thus,~$X$
may be given by two equations of the form
\begin{eqnarray*}
-c_0 T_{i_0}^2 +\, c_1 T_{i_1}^2 +\, c_2 T_{i_2}^2 +\, c_3 T_{i_3}^2 & = & 0 \, , \\
-d_0 T_{j_0}^2 + d_1 T_{j_1}^2 + d_2 T_{j_2}^2 + d_3 T_{j_3}^2 & = & 0 \, ,
\end{eqnarray*}
for
$c_k, d_k > 0$,
$i_0 \neq j_0$,
and~$\{i_0,\ldots,i_3\} \cup \{j_0,\ldots,j_3\} = \{0,\ldots,4\}$.
These~equations imply
$x_{i_0} \neq 0$
and~$x_{j_0} \neq 0$
for every real point
$x \in X(\bbR)$.
In~particular,
$X(\bbR)$
has at least two connected components, given by the two possible signs
of~$x_{i_0}/x_{j_0}$.
Clearly,~these two components are interchanged under the operation
of~$\Aut'(X)$.

We~finally note that a real degree four del Pezzo surface cannot have more than two connected components~\cite[Chapter~III, Theorem~3.3]{Silh}.
}
\eop
\end{lem}

\begin{rem}
The~stronger statement that if
$X(\bbR)$
splits into two connected components then the operation
of\/~$\Aut'(X)$
interchanges them is true, as~well.

Indeed,~by blowing up a real point not lying on any exceptional curve, one obtains a real cubic surface that has two connected components. According to L.~Schl\"af\/li~\cite[pp.~114f.]{Sch}, there are exactly five real types of real cubic surfaces and those correspond in modern language to the four conjugacy classes of
order-$2$
subgroups
in~$W(E_6)$
together with the trivial~group. Only~for one of these five cases, the Brauer group is nontrivial \mbox{\cite[Appendix, Table~2]{Ja}}, it is isomorphic
to~$(\bbZ/2\bbZ)^2$
then, and that is the single case in which the surface is~disconnected. We~will not make use of this~observation.
\end{rem}

\noindent
{\bf Proof of Theorem~\ref{theo2}.}
Let
$X$
be a diagonalizable degree four del Pezzo
surface over~$\bbQ$
that has an adelic point and a Brauer class
$\alpha \in \Br(X)$
working at the infinite~place. We~note that, since
$X$
has an adelic point, it clearly has a real~point. The~local evaluation
$\ev_{\alpha,\infty}(x)$
for~$x \in X(\bbR)$
is defined using the restriction homomorphism
$x^*\colon \Br(X) \rightarrow \Br(\bbR)$,
which factors via
$\Br(X_\bbR)$.
Hence,~non-constancy of
$\ev_{\alpha,\infty}$
implies that the restriction
$\alpha_\bbR \in \Br(X_\bbR)/\Br(\bbR)$
is a nonzero~class.

In~this case, Lemma~\ref{XR} shows that
$X(\bbR)$
splits into two connected components. Moreover, there exists an automorphism
$\sigma \in \Aut'(X_\bbR) = \Aut'(X)$
interchanging these.
Since~$\ev_{\alpha,\infty}$
is locally constant, this implies that
$\sigma$
is a witness for the non-constancy of the local evaluation map
at~$\infty$.
In~other words, the natural homomorphism
$i_\alpha\colon \Aut'(X) \rightarrow \Br(\bbQ)$
has in its image a class
$i_\alpha(\sigma)$
with a non-zero component at~infinity.

According to global class field theory~\cite[Section~10, Theorem~B]{Ta},
$i_\alpha(\sigma)$
necessarily has a nonzero component at a second place
$\nu \neq \infty$.
Consequently,
$\alpha$
works at the
place~$\nu$,
too, which implies the~claim.
\eop

\section{Surfaces with a Brauer class working at \\ a prescribed set of places}

The goal of this section is to prove Theorem~\ref{theo3}. We~distinguish between the cases $\#S > 1$,
$\#S = 1$,
and~$S = \emptyset$.
The~family below will serve us in all~cases.

\subsection{A family of degree four del Pezzo surfaces}

For~$D, A_1, A_2, B \in \bbQ$,
let
$S := S^{(D;A_1,A_2,B)} \subset \Pb^4_\bbQ$
be given by the system of~equations
\begin{eqnarray}
\label{diagf1}
 -A_1(T_0-T_1)(T_0+T_1) & = & T_3^2 - DT_4^2 \, , \\
\label{diagf2}
 -A_2(T_0-T_2)(T_0+T_2) & = & T_3^2 - B^2DT_4^2 \, .
\end{eqnarray}

\begin{theo}
\label{fam}
Let~$D, A_1, A_2, B$
be nonzero rational numbers.

\begin{abc}
\item[{\rm A.a) }]\addtocounter{abc}{1}
Then\/~$S$
is not a~cone. The
degree-$5$
scheme\/
$\calS_X$
has a point at infinity and four others, which are the roots of a completely reducible polynomial of degree four having discriminant\/
$\Delta := A_1^2 (A_1-A_2)^2 (A_1B^2-A_2)^2 B^4 (B-1)^2 (B+1)^2/A_2^6 B^{12}$.
\item
$S$
has the\/
\mbox{$\bbQ$-rational}
point\/~$(1\!:\!1\!:\!1\!:\!0\!:\!0) \in X(\bbQ)$.
\item
If\/~$\Delta \neq 0$
then the five
rank-$4$
discriminants are, up to perfect square factors, given by
$D$,
$D$,
as well as\/
$DA_1A_2(A_1-A_2)(B^2-1)$,
$A_1A_2(A_1B^2-A_2)(B^2-1)$,
and\/
$D(A_1-A_2)(A_1B^2-A_2)$.
\item[{\rm B.a) }]\addtocounter{abc}{-2}
There~is a Brauer class\/
$\alpha \in \Br(X)$
extending that of the quaternion algebra\/
$\smash{\big( \bbQ(\sqrt{D})(X),\tau,\frac{T_0+T_1}{T_0+T_2} \big)}$
over the function
field\/~$\bbQ(X)$.
\item
Moreover,~one has\/
$\ev_{\alpha,\nu}(x) = 0$
for\/
$x = (1\!:\!1\!:\!1\!:\!0\!:\!0)$
and
every\/~$\nu \in \Omega$.
\item
At~an arbitrary place\/
$\nu \in \Omega$,
the local evaluation map\/
$\ev_{\alpha,\nu}$
is constant if one of the following conditions~holds.
\item[ $\bullet$ ]\looseness-1
$\nu = p$
is a finite place,
$p \neq 2$,
and\/
$p$
divides neither\/
$D$,
nor\/~$A_1$,
nor\/~$A_2$,
nor\/~$B$.
\item[ $\bullet$ ]
$\nu = p$~splits
in\/~$\smash{\bbQ(\sqrt{D})}$,
or\/
$\nu = \infty$
and\/~$D > 0$.
\item[ $\bullet$ ]
$D$~is
square-free,
$\nu = p$
is a finite place,
$p|D$,
$p \neq 2$,
$\gcd(B,D) = 1$,
$\smash{(\frac{-A_1}p) = 1}$,
and\/
$A_1 \equiv A_2 \pmodulo p$.
\item
At~a place\/
$\nu$,
the local evaluation map\/
$\ev_{\alpha,\nu}$
cannot be constant if\/
$(-A_1,D)_\nu = -1$
or\/~$(-A_2,D)_\nu = -1$.
\end{abc}\medskip

\noindent
{\bf Proof.}
{\em
A.a) and~c) are standard calculations, while b)~is directly~checked. Moreover,~B.a) is a direct application of Proposition~\ref{Brauerklasse}.a) and the assertion of~b) follows from the fact that
$\frac{x_0+x_1}{x_0+x_2} = 1$
for~$x = (1\!:\!1\!:\!1\!:\!0\!:\!0)$.\smallskip

\noindent
B.c)
The~sufficiency of the first condition is Criterion~\ref{VV}, while that of the second was shown in Proposition~\ref{Brauerklasse}.b). In~order to establish the sufficiency of the third, we argue as~follows.

First~of all, the prime
$p$
ramifies
in~$\smash{\bbQ(\sqrt{D})}$.
A~\mbox{$p$-adic}
unit
$u \in \bbQ_p$
is a local norm from
$\smash{\bbQ(\sqrt{D})}$
if and only
if~$(u \bmod p) \in \bbF_{\!p}^*$
is a~square. Moreover,~we note that
$\smash{(\frac{-A_1}p) = (\frac{-A_2}p) = 1}$
implies that
each of the four rational functions
$\smash{\frac{T_0 \pm T_1}{T_0 \pm T_2}}$
defines the Brauer
class~$\alpha$.

Let~now
$x \in X(\bbQ_p)$
be any
\mbox{$p$-adic}
point. Normalize~the coordinates
$x_0,\ldots,x_4$
such that each is a
\mbox{$p$-adic}
integer and at least one is a~unit.
If~$p|x_0$
and
$p|x_1$
or
$p|x_0$
and
$p|x_2$
then the equations imply that all coordinates must be divisible
by~$p$,
a~contradiction.
Modulo~$p$,
we have
$\overline{x}_0^2 - \overline{x}_1^2 = \overline{x}_0^2 - \overline{x}_2^2$,
hence~$\overline{x}_1 = \pm \overline{x}_2$,
which implies that one of the four quotients
$\smash{\frac{x_0 \pm x_1}{x_0 \pm x_2}}$
is congruent
to~$1$
modulo~$p$,
and therefore a~norm.\smallskip

\noindent
B.d)
We~note first that
$X$
has
\mbox{$X(\bbQ_\nu)$-rational}
points such that
$x_0 \neq \pm x_1$
and
$x_0 \neq x_2$.
Indeed,~putting
$x_0 := 1$
and choosing
$x_3$
and
$x_4$
sufficiently close
to~$0$
in the
\mbox{$\nu$-adic}
topology, we see that (\ref{diagf1}) and (\ref{diagf2}) become soluble when viewed as equations for
$x_1$
and
$x_2$,~respectively.

Now,~let us assume without loss of generality that
$(-A_1,D)_\nu = -1$.
Then~the automorphism
$\sigma\colon (T_0:\ldots:T_4) \mapsto (T_0\!:\!(-T_1)\!:\!T_2\!:\!T_3\!:\!T_4)$
changes the rational function
$\smash{\frac{T_0+T_1}{T_0+T_2}}$
by a factor of
$$\textstyle \frac{T_0-T_1}{T_0+T_1} = -\frac1{A_1} \frac{T_3^2-DT_4^2}{(T_0+T_1)^2} \, ,$$
which is a
\mbox{$\nu$-adic}
non-norm
from~$\smash{\bbQ(\sqrt{D})}$
since
$(-A_1,D)_\nu = -1$.
This~shows that
$i_\alpha(\sigma)$
has a nonzero component at
$\nu$,
i.e.\ that the automorphism
$\sigma$
witnesses the non-constancy of the local evaluation
map~$\ev_{\alpha,\nu}$.
}
\eop
\end{theo}

\subsection{More than one place}

Let~$S \subset \Omega$
be a given set that consists of at least two~places. We~incorporate the notation
$\{p_1,\ldots,p_r\} = S \!\setminus\! \{2,\infty\}$.

To~construct a diagonalizable degree four del Pezzo surface such that a nontrivial Brauer class works exactly at the places
in~$S$,
we first choose a square-free integer
$D \neq 0$
satisfying the following~conditions.

\begin{iii}
\item[ $\bullet$ ]
$D > 0$
if and only if
$\infty \not\in S$.
\item[ $\bullet$ ]
$D \equiv 3 \pmodulo 4$
in the case that
$2 \in S$,
and
$D \equiv 1 \pmodulo 8$
when
$2 \not\in S$.
\item[ $\bullet$ ]
$D$
is divisible by
$p_1,\ldots,p_r$
and has exactly one further prime divisor, which we
call~$q$.
\end{iii}

\noindent
That~such a choice
of~$D$
is possible follows immediately from the fact that there are infinitely many primes in every odd residue class
modulo~$8$.

Now~write
$S = S_1 \cup S_2$
as a union of two not necessarily disjoint subsets of even~size. This~is possible, because of
$\#S \geq 2$.
In~addition, we may put
$2$
into both subsets in case it occurs as an element
of~$S$,
and the same
for~$\infty$.

Next,~we choose primes
$A_1 \neq A_2$
not
dividing~$D$
such that, for
$i = 1,2$,
\begin{equation}
\label{Hilbcond}
(-A_i,D)_\nu = -1 \quad \Longleftrightarrow \quad \nu \in S_i \, .
\end{equation}
To~see that this may be achieved, we observe at first that
$(-A_i,D)_\nu = 1$
for all places
$\nu \neq 2,\infty; p_1,\ldots,p_r,q$,
and~$A_i$.
The~requirement
at~$\nu=2$
may be realized by choosing
$A_i \equiv 1 \pmodulo 4$,
the~condition
at~$\nu=\infty$
is implied by the choice that
$A_i$
is~positive. Furthermore,~we require
$$\textstyle \big( \frac{-A_i}{p_j} \big) = \left\{
\begin{array}{rl}
-1 & {\rm if~} p_j \in S_i \,, \\
 1 & {\rm otherwise} \,,
\end{array}
\right.$$
and
$\smash{(\frac{-A_i}{q}) = 1}$.
Let~us impose, in addition, the condition that
\begin{equation}
\label{extra}
A_2 \equiv A_1 \pmodulo q \, .\smallskip
\end{equation}

All~these are congruence conditions modulo distinct odd primes. Therefore,~the existence of a prime
$A_1$
satisfying~(\ref{Hilbcond}) for all places except, possibly,
$A_1$~itself,
is implied by Dirichlet's Theorem on primes in arithmetic~progressions. Moreover,~as
$\#S_i$
is even, we have
$(-A_i,D)_{A_i} = 1$
by the Hilbert reciprocity law~\cite[Chapter~VI, Theorem~8.1]{Ne}.

In~a completely analogous manner, Dirichlet's Theorem and the Hilbert reciprocity law imply the existence of a prime
$A_2 \neq A_1$
fulfilling (\ref{Hilbcond}) and~(\ref{extra}).\medskip

We~may now formulate the main result of this~paragraph.

\begin{theo}
\label{ge2}
Let~the integers\/
$D$,
$A_1$,
and\/
$A_2$
be chosen as~above.

\begin{abc}
\item
Then, for every
integer\/~$B \geq 2$,
the surface\/
$X \subset \Pb^4_\bbQ$
given~by
\begin{eqnarray*}
 -A_1(T_0-T_1)(T_0+T_1) & = & T_3^2 - DT_4^2 \, , \\
 -A_2(T_0-T_2)(T_0+T_2) & = & T_3^2 - B^2DT_4^2
\end{eqnarray*}
is nonsingular and has a\/
\mbox{$\bbQ$-rational}
point.
\item
There~is a Brauer class\/
$\alpha \in \Br(X)$
extending that of the quaternion algebra\/
$\smash{\big( \bbQ(\sqrt{D})(X),\tau,\frac{T_0+T_1}{T_0+T_2} \big)}$
over the function
field\/~$\bbQ(X)$.
\item
The~Brauer class\/
$\alpha$
works at every place
$\nu \in S$.
If\/~$B$
is a prime number that splits
in\/~$\smash{\bbQ(\sqrt{D})}$
then
$\alpha$
does not work at any~other place.
\item
There~are infinitely many prime
numbers\/~$B$
splitting
in\/~$\smash{\bbQ(\sqrt{D})}$
such that\/
$\Br(X)/\Br(\bbQ) \cong \bbZ/2\bbZ$.
\end{abc}\medskip

\noindent
{\bf Proof.}
{\em
a)~follows from Theorem~\ref{fam}.A.a) and~b). b)~is Theorem~\ref{fam}.B.a).\smallskip

\noindent
c)
Our choices of
$A_1$,
$A_2$,
and~$D$
guarantee that Theorem~\ref{fam}.B.d) applies to every
$\nu \in S$.
On~the other hand, as
$B$
is a prime that splits in 
$\smash{\bbQ(\sqrt{D})}$,
Theorem~\ref{fam}.B.c) shows constancy of the evaluation map
$\ev_{\alpha,\nu}$
for all other~places.\smallskip

\noindent
d)
According~to Fact~\ref{order4}, in order to exclude the option
$\Br(X)/\Br(\bbQ) \cong (\bbZ/2\bbZ)^2$,
we have to choose the parameter
$B$
such that neither of the terms
$$A_1A_2(A_1-A_2)(B^2\!-1), \; DA_1A_2(A_1B^2\!-A_2)(B^2\!-1), \; {\rm and} \; (A_1-A_2)(A_1B^2\!-A_2)$$
is a perfect~square. By~Siegel's Theorem on integral points on elliptic curves~\cite[Theorem IX.4.3]{Silv}, the term in the middle is a square only finitely many~times. The~two others lead to Pell-like equations, the integral solutions of which are known to have exponential~growth, cf.\ for example \cite[Chapter~XXXIII, \S\S15--18]{Ch}. The~assertion~follows.
}
\eop
\end{theo}

\subsection{No place}

It~is not at all hard to write down a diagonalizable degree four del Pezzo
surface~$X$
such that a nontrivial Brauer class works at no place at~all.

\begin{ex}
Let\/~$X \subset \Pb^4_\bbQ$
be the surface given~by
\begin{eqnarray*}
    -(T_0-T_1)(T_0+T_1) & = & T_3^2 - 17T_4^2 \, , \\
 -103(T_0-T_2)(T_0+T_2) & = & T_3^2 - 68T_4^2 \, .
\end{eqnarray*}
Then~$X$
is nonsingular
and~$X(\bbQ) \neq \emptyset$.
Moreover,
$\Br(X)/\Br(\bbQ) \cong \bbZ/2\bbZ$
but the nontrivial Brauer class works at no~place.\medskip

\noindent
{\bf Proof.}
The~first two assertions follow from Theorem~\ref{fam}.A.b) and~a). The~discriminants of the five
rank-$4$
forms are, up to square factors,
$17$,
$17$,
$66$,
$206$,
and~$3399$
such that, by Facts~\ref{order2}.a) and~\ref{order4}, we have
$\Br(X)/\Br(\bbQ) \cong \bbZ/2\bbZ$.

Let~$\alpha \in \Br(X)$
be~nontrivial. By~Theorem~\ref{fam}.B.c), the local evaluation map is constant at all places
$\nu \neq 2,17,103$,
and~$\infty$.
Moreover,~it is constant
at~$\nu=\infty$
as the field
$\smash{\bbQ(\sqrt{17})}$
is real-quadratic. Constancy at
$\nu=2$
and~$103$
is clear, too, since these primes split
in~$\smash{\bbQ(\sqrt{17})}$.
Finally,~$\ev_{\alpha,17}$
is constant as
$\smash{(\frac{-1}{17}) = 1}$
and
$103 \equiv 1 \pmodulo {17}$.
\eop
\end{ex}

\subsection{Exactly one place}

The~examples here are necessarily a bit different, as the 16 automorphisms must not witness the non-constancy of the evaluation~map. We~may nonetheless work with the family from Theorem~\ref{fam}.

\begin{ex}
Let~$l$
be a prime number such that
$l \equiv 3 \pmodulo 4$.
Choose~a prime
$D \equiv 1 \pmodulo 8$
such that
$\smash{(\frac{D}l) = -1}$
and another prime
$A > l$
such that
$A \equiv 1 \pmodulo D$
and
$(A^2-1)(A^2-l^2)$
is a non-square.\smallskip

\noindent
Then~the surface
$X \subset \Pb^4_\bbQ$
given~by
\begin{eqnarray*}
    -(T_0-T_1)(T_0+T_1) & = & T_3^2 - DT_4^2 \, , \\
 -A^2(T_0-T_2)(T_0+T_2) & = & T_3^2 - l^2DT_4^2
\end{eqnarray*}
is nonsingular and has a
\mbox{$\bbQ$-rational}
point.
Moreover,~$\Br(X)/\Br(\bbQ) \cong \bbZ/2\bbZ$
and the nontrivial class works exactly at the
place~$l$.\medskip

\noindent
{\bf Proof.}
We~first note that the restrictions on
$D$
and~$A$
are easy to fulfill due to Dirichlet's and Siegel's~theorems.
Furthermore,~the first three assertions follow directly from Theorem~\ref{fam}.A.a) and~b), as well as Facts~\ref{order2}.a) and \ref{order4}.
The~nontrivial Brauer class
$\alpha \in \Br(X)$
may be understood as an extension of the quaternion algebra
$\smash{\big( \bbQ(\sqrt{D})(X),\tau,\frac{T_0+T_1}{T_0+T_2} \big)}$
over
$\bbQ(X)$
to the whole
scheme~$X$.
Moreover,~any of the quotients
$\smash{\frac{T_0 \pm T_1}{T_0 \pm T_2}}$
defines the same Brauer~class.
Theorem~\ref{fam}.B.c) implies that the local evaluation map is constant at all places
$\nu \neq l$.\smallskip

\noindent
{\em Non-constancy of\/
$\ev_{\alpha,l}$:\/}
Note that
$l$
is an inert prime,
since~$\smash{(\frac{D}l) = -1}$.
An~element
$\smash{u \in \bbQ_l^*}$
is a local norm from
$\smash{\bbQ(\sqrt{D})}$
if and only
if~$\nu_l(u)$
is~even.

For~$\underline{x} = (1\!:\!1\!:\!1\!:\!0\!:\!0)$,
we have
$\ev_{\alpha,l}(\underline{x}) = 0$
by~Theorem~\ref{fam}.B.b). On~the other hand, the substitutions
$T_0 = lT'_0$,
$T_1 = T'_1$,
$T_2 = lT'_2$,
$T_3 = lT'_3$,
and
$T_4 = T'_4$
yield a different model
$X'$
of~$X$
that is given~by
\begin{eqnarray*}
 -(lT'_0-T'_1)(lT'_0+T'_1) & = & l^2T'_3 {}^2 - DT'_4 {}^2 \, , \\
 -A^2(\phantom{l}T'_0-T'_2)(\phantom{l}T'_0+T'_2) & = & \phantom{l^2} T'_3 {}^2 - DT'_4 {}^2 \, .
\end{eqnarray*}
Moreover,
$\smash{\frac{T_0+T_1}{T_0+T_2} = \frac{lT'_0+T'_1}{lT'_0+lT'_2} = \frac1l \frac{lT'_0+T'_1}{T'_0+T'_2}}$.
It~suffices to find a
\mbox{$\bbQ_l$-rational}
point
on~$X'$~so
that
$\smash{\frac{lT'_0+T'_1}{T'_0+T'_2}}$
is a
\mbox{$l$-adic}~unit.

The~reduction
of~$X'$
modulo
$l$~is
given~by
\begin{eqnarray*}
                             T'_1 {}^2 & = & {\phantom{T'_3 {}^2}} \,\,-\!\!\overline{D}T'_4 {}^2 \, , \\
 -\overline{A}^2(T'_0-T'_2)(T'_0+T'_2) & = & T'_3 {}^2 - \overline{D}T'_4 {}^2 \, .
\end{eqnarray*}
We~observe that the first equation has a nontrivial solution, as
$\smash{(\frac{D}l) = -1}$
and
$l \equiv 3 \pmodulo 4$
together imply that 
$(-\overline{D}) \in \bbF_{\!l}^*$
is a~square.
Let~$\rho \in \bbF_{\!l}^*$
be one of its square~roots.

The~Jacobian matrix associated to a point
$x \in X'(\bbF_{\!l})$~is
$$
\left(
\begin{array}{ccccc}
0 & 2x_1 & 0 & 0 & 2\overline{D}x_4 \\
-2\overline{A}^2\!x_0 & 0 & 2\overline{A}^2\!x_2 & -2x_3 & 2\overline{D}x_4
\end{array}
\right) ,
$$
which shows that points such that
$x_1 \neq 0$
are non-singular.
For~instance, there is the non-singular
\mbox{$\bbF_{\!l}$-rational}
point
$$\textstyle x = (\frac{\overline{A}^2+\overline{D}}{2\overline{A}^2}:\rho:\frac{\overline{A}^2-\overline{D}}{2\overline{A}^2}:0:1) \in X'(\bbF_{\!l}) \, .$$
Here,
$\smash{\frac{x_1}{x_0+x_2} = \rho}$
is a nonzero element
in~$\bbF_{\!l}$.
Hence,~for every
\mbox{$l$-adic}
point that
lifts~$x$,
the value of the quotient
$\smash{\frac{lT'_0+T'_1}{T'_0+T'_2}}$
is an
\mbox{$l$-adic}
unit as required. The~assertion follows.%
\eop
\end{ex}

In order to provide the corresponding example in the
$l \equiv 1 \pmodulo 4$
case, we need the following~lemma.

\begin{lem}
\label{conicpt}
Let\/~$\bbF_{\!l}$
be a finite field of
characteristic\/~$\neq\! 2$.
Then~there exists an element\/
$\sigma \in \bbF_{\!l}$
such that\/
$2(1+\sigma^2)$
is a non-square
in\/~$\bbF_{\!l}$.\medskip

\noindent
{\bf Proof.}
{\em
If~$l=3$
then
put~$\sigma := 0$.
Otherwise,~i.e.~for
$l \geq 5$,
let~$c \in \bbF_{\!l}$
be any non-square. The~equation
$cT_0^2 = 2(T_1^2 + T_2^2)$
defines a conic over
$\bbF_{\!l}$,
which has exactly
$(l+1)$
\mbox{$\bbF_{\!l}$-rational}
points. Among~them, at most four have
$x_1 = 0$
or~$x_0 = 0$.
For~the others,
$\sigma := x_2/x_1$
fulfills the required~condition.
}
\eop
\end{lem}

\begin{ex}
Let~$l$
be a prime number that is either
$l \equiv 1 \pmodulo 4$
or~$l=2$.
If\/~$l=2$
choose\/~$B := 2$,
otherwise let\/
$B$
be an odd prime number that is split
in\/~$\smash{\bbQ(\sqrt{\mathstrut l})}$
and such that neither
$l(l-1)(B^2-1)$,
nor~$(lB^2-1)(B^2-1)$,
nor~$(l-1)(lB^2-1)$
is a perfect~square.\smallskip

\noindent
Then~the surface
$X \subset \Pb^4_\bbQ$
given~by
\begin{eqnarray}
\label{exfl1}
-l(T_0-T_1)(T_0+T_1) & = & T_3^2 - lT_4^2 \, , \\
\label{exfl2}
 -(T_0-T_2)(T_0+T_2) & = & T_3^2 - B^2lT_4^2
\end{eqnarray}
is nonsingular and has a
\mbox{$\bbQ$-rational}
point.
Moreover,~$\Br(X)/\Br(\bbQ) \cong \bbZ/2\bbZ$
and the nontrivial class works exactly at the
place~$l$.\medskip

\noindent
{\bf Proof.}
Once~again, the restrictions on
$B$
are easy to fulfill due to Dirichlet's and Siegel's~theorems.
Furthermore,~the first three assertions follow directly from Theorem~\ref{fam}.A.a) and~b) as well as Facts~\ref{order2}.a) and Fact~\ref{order4}.
There is a Brauer class
$\alpha \in \Br(X)$,
which may be understood as an extension of the quaternion algebra
$\smash{\big( \bbQ(\sqrt{l})(X),\tau,\frac{T_0+T_1}{T_0+T_2} \big)}$
over
$\bbQ(X)$
to the whole
scheme~$X$.
Moreover,~any of the quotients
$\smash{\frac{T_0 \pm T_1}{T_0 \pm T_2}}$
defines the same Brauer~class.
Finally,~Theorem~\ref{fam}.B.c) implies that the local evaluation map is constant at all places
$\nu \neq 2,l$.
Thus,~all which remains to be shown is that
$\ev_{\alpha,l}$
is non-constant and that
$\ev_{\alpha,2}$
is constant in the case
$l \equiv 1 \pmodulo 4$.\smallskip

\noindent
{\em Non-constancy of\/
$\ev_{\alpha,2}$
for\/~$l=2$:\/}
For~$\smash{\bbQ(\sqrt{2})}$,
the prime
$p=2$
is ramified. A
\mbox{$2$-adic}
unit~$u$
is a local norm from
$\smash{\bbQ(\sqrt{2})}$
if and only
if~$u \equiv \pm1 \pmodulo 8$.

For~$\underline{x} = (1\!:\!1\!:\!1\!:\!0\!:\!0)$,
we have
$\ev_{\alpha,2}(\underline{x}) = 0$
by~Theorem~\ref{fam}.B.b). On~the other hand, there is the
\mbox{$2$-adic}
point
$\smash{x = (1\!:\!0\!:\!\sqrt{-7}\!:\!0\!:\!1) \in X(\bbQ_2)}$.
Observe~that
$(-7) \equiv 1 \pmodulo 8$
implies that
$(-7)$
is a square
in~$\bbQ_2$.
Moreover,~we may choose
$\sqrt{-7} \in 5 + 16\bbZ_2$
since
$5^2 \equiv -7 \pmodulo {32}$.
Then~$\smash{\frac{x_0+x_1}{x_0+x_2} = 1/(1+\sqrt{-7})}$,
which is in the residue class
$\frac12 \cdot (3 \bmod 8)$.
Consequently,
$\ev_{\alpha,2}(x) = \frac12$.\smallskip

\noindent
{\em Non-constancy of\/
$\ev_{\alpha,l}$
for\/~$l \equiv 1 \pmodulo 4$:\/}
For~$\smash{\bbQ(\sqrt{\mathstrut l})}$,
the prime
$l$
is ramified.
A~\mbox{$l$-adic}
unit~$u$
is a local norm from
$\smash{\bbQ(\sqrt{\mathstrut l})}$
if and only
if~$(u \bmod l) \in \bbF_{\!l}^*$
is a~square.

For~$\underline{x} = (1\!:\!1\!:\!1\!:\!0\!:\!0)$,
we have
$\ev_{\alpha,l}(\underline{x}) = 0$.
On~the other hand, the substitution
$T_3 = lT'_3$
yields a different model
$X'$
of~$X$
that is given~by
\begin{eqnarray*}
 (T_0-T_1)(T_0+T_1) & = & T_4^2 - lT'_3 {}^2 \, , \\
 (T_0-T_2)(T_0+T_2) & = & l(B^2T_4^2 - lT'_3 {}^2) \, .
\end{eqnarray*}
The~reduction
of~$X'$
modulo
$l$~is
given~by
\begin{eqnarray*}
 (T_0-T_1)(T_0+T_1) & = & T_4^2 \, , \\
 (T_0-T_2)(T_0+T_2) & = & 0 \, .
\end{eqnarray*}
From~this, we see that the Jacobian matrix associated to a point
$x \in X'(\bbF_{\!l})$
is
$$
\left(
\begin{array}{ccccc}
 2x_0 & -2x_1 & 0 & 0 & -2x_4 \\
 2x_0 & 0 & -2x_2 & 0 & 0
\end{array}
\right) ,
$$
which shows that points such that
$x_0 = x_2 \neq 0$
are non-singular.
For~instance, for any
$\sigma \in \bbF_{\!l}$
such that
$\sigma^2 \neq -1$,
there is the non-singular
\mbox{$\bbF_{\!l}$-rational}
point
$$\textstyle x = \big( (1+\sigma^2)\!:\!2\sigma\!:\!(1+\sigma^2)\!:\!0\!:\!(1-\sigma^2) \big) \in X(\bbF_{\!l}) \, .$$
If,~moreover,
$\sigma$
is chosen as in Lemma~\ref{conicpt} then
$\smash{\frac{x_0+x_1}{x_0+x_2} = \frac{(1+\sigma)^2}{2(1+\sigma^2)}}$,
which is a non-square. Then,~for every
\mbox{$l$-adic}
point that
lifts~$x$,
the local evaluation map has
value~$\frac12$.\smallskip

\noindent
{\em Constancy of\/
$\smash{\ev_{\alpha,2}}$
for\/~$l \equiv 1 \pmodulo 4$:\/}
If~$l \equiv 1 \pmodulo 8$
then the prime
$p=2$
is split in
$\smash{\bbQ(\sqrt{\mathstrut l})}$
and there is nothing to~prove.

On~the other hand, assume that
$l \equiv 5 \pmodulo 8$,
in which case
$2$~is
an inert~prime. Then~an element
$u \in \bbQ_2^*$
is a local norm from
$\smash{\bbQ(\sqrt{\mathstrut l})}$
if and only
if~$\nu_2(u)$
is~even. Moreover,~any
\mbox{$2$-adic}
point
$x \in X(\bbQ_2)$
may be represented by coordinates
$x_0,\ldots,x_4$
that are
\mbox{$2$-adic}
integers, at least one of which is a~unit.
It~is now a routine matter to determine all quintuples of residues
modulo~$8$,
one of which is odd, that satisfy the system (\ref{exfl1},\ref{exfl2}) of equations
modulo~$8$.
From~the list obtained, one already sees that
$\ev_{\alpha,2}(x) = 0$
in each~case. We~leave the details to the~reader.
\eop
\end{ex}

\section{A consequence concerning del Pezzo surfaces of low degree}

The aim of this section is to prove Theorem~\ref{theo4}. For~this, we need some technical facts about point sets
on~$\Pb^2$
that lie in general position, which should be well-known among experts. We decided to include them due to the lack of a suitable~reference.

Let~$K$
be an algebraically closed~field. One~says that
$r \leq 8$
distinct points
$P_1,\ldots,P_r \in \Pb^2(K)$
lie {\em in general position,} if no three of them lie an a line, no six on a conic, and no eight on a cubic that has a singular point at one of~them.
It~is well known, cf.~\cite[Proposition 8.1.25]{Do}, \cite[Section~3]{Is}, or \cite[Theorem~1]{De}, that the blow-up of
$\Pb^2_K$
in
$P_1,\ldots,P_r$
is a del Pezzo surface if and only if these points lie in general position.

\begin{lem}
\label{cub7}
Let\/~$K$
be an algebraically closed field and\/
$P_1,\ldots,P_7 \in \Pb^2(K)$
be seven points in general~position. Then~there exists a nonsingular cubic through all of~them.\medskip

\noindent
{\bf Proof.}
{\em
By~\cite[Corollary~V.4.4.a) and Proposition~V.4.3]{Ha}, the linear system of all cubics through
$P_1,\ldots,P_7$
is two-dimensional and has no unassigned base~points. Blowing~up
$P_1,\ldots,P_7$,
we find a two-dimensional linear system on a non-singular surface that is base point~free.

A~version of Bertini's theorem~\cite[Corollary III.10.9]{Ha} shows that the generic element of this linear system is~smooth. Projecting down
to~$\Pb^2$,
we see that no singular points may occur, except for
$P_1,\ldots,P_7$.

However,~according to \cite[Theorem~26.3]{Man}, there are only seven cubics through
$P_1,\ldots,P_7$
that are singular at one of~them.
Therefore,~the generic member must be nonsingular, as~required.
}
\eop
\end{lem}

\begin{prop}
\label{genpos}
Let\/~$K$
be an algebraically closed field and\/
$P_1,\ldots,P_5 \in \Pb^2(K)$
be five points in general~position. Then~there exist\/
$P_6,P_7,P_8 \in \Pb^2(K)$
such that\/
$P_1,\ldots,P_8 \in \Pb^2(K)$
lie in general~position.\medskip

\noindent
{\bf Proof.}
{\em First~choose
$P_6 \in \Pb^2(K)$
not lying on any of the 10~lines through two of the five points, and neither on the conic through all of~them.
Next,~choose
$P_7 \in \Pb^2(K)$
outside the 15~lines through two of the six points and outside the six conics through five of~them.
All~these requirements exclude only a one-dimensional subset
of~$\Pb^2$.

When~we finally choose
$P_8$,
we have to be slightly more careful, as singular cubics need to be taken into consideration. Clearly,~we have to choose
$P_8 \in \Pb^2(K)$
outside the 21~lines through two of the seven points, outside the 21~conics through five of them, and outside the seven cubics through the points
$P_1,\ldots,P_7$
that are singular at one of~them. This~excludes, once again, only a one-dimensional subset
of~$\Pb^2$.

There~is one final condition.
$P_8$
must not be the singular point of a cubic through
$P_1,\ldots,P_7$.
To~analyze this requirement, recall~\cite[Corollary~V.4.4.a)]{Ha} that the linear system of all cubics through these seven points is two dimensional. By~Lemma~\ref{cub7}, it contains a nonsingular curve. Thus,~singular curves form an at most one-dimensional subfamily. Moreover,~none of these cubics may have multiple components, so that the total set of singular points occurring is at most one-dimensional. This~completes the~proof.
}
\eop
\end{prop}

\begin{coro}
\label{blow}
Let\/~$K$
be an algebraically closed field and\/
$X$
be a del Pezzo surface of degree four
over\/~$K$.
Then~the set of all\/
$(P,Q\!,R) \in X^3(K)$
such that\/
$\Bl_{\{P,Q,R\}}(X)$
is a del Pezzo surface of degree one is Zariski open and non-empty.\medskip

\noindent
{\bf Proof.}
{\em
$X$
is isomorphic to
$\Pb^2$,
blown up in five points in general~position. Proposition~\ref{genpos} shows that the set considered is non-empty. Moreover,~being in general position is an open condition, so that Zariski openness is clear.
}
\eop
\end{coro}\medskip

\noindent
{\bf Proof of Theorem~\ref{theo4}.}
For~$d=4$,
this is Theorem~\ref{theo1}. It~yields a degree four del Pezzo surface
$X$
having a
\mbox{$\bbQ$-rational}
point such that
$\Br(X)/\Br(\bbQ) \cong \bbZ/2\bbZ$
and the nontrivial Brauer class works exactly at the places
in~$S$.
We~note that, by~\cite[Theorem~0.1]{SSk}, the
\mbox{$\bbQ$-rational}
points
on~$X$
are Zariski~dense.

Moreover,~Brauer groups do not change under blowup and the local evaluation maps are compatible in the sense that
$\ev_{\alpha,p}(\pi(x)) = \ev_{\pi^*\!\alpha,p}(x)$.
Thus,~given some integer
$1 \leq d < 4$,
let us blow up
$X$
in
$(4-d)$
\mbox{$\bbQ$-rational}
points. This~clearly yields a
surface~$Y$
over~$\bbQ$
that has a
\mbox{$\bbQ$-rational}
point and fulfills the conditions that
$\Br(Y)/\Br(\bbQ) \cong \bbZ/2\bbZ$
and the nontrivial Brauer class works exactly at the places
in~$S$.

As~$\smash{K_Y^2 = 4 - (4-d) = d}$,
it only remains to ensure that we may choose the blow-up points in such a way that
$Y$
becomes a del Pezzo surface. In~view of~\cite[Corollary~24.5.2.i)]{Man}, it suffices to do this in the case
when~$d=1$.

For~this, let us fix an algebraic closure
$\overline\bbQ$
and view
$\bbQ$
as a subfield of~it. By~Corollary~\ref{blow}, we know that the set of all
$(P,Q\!,R) \in X^3(\overline\bbQ)$
that yield a del Pezzo surface is Zariski open and non-empty. On~the other hand, since
$X(\bbQ)$
is Zariski dense
in~$X(\overline\bbQ)$,
we also have that
$X^3(\bbQ)$
is Zariski dense
in~$X^3(\overline\bbQ)$.
As~a non-empty open subset and a dense one necessarily have a point in common, the proof is complete.\eop

\setlength\parindent{0mm}

\frenchspacing
\begin{thebibliography}{BMS}
\bibitem[AG]{AG}
Auslander, M.\ and Goldman, O.: The Brauer group of a commutative ring, {\em Trans.\ Amer.\ Math.\ Soc.} {\bf 97}\br(1960)\brr367--409
\bibitem[BSD]{BSD}
Birch, B.\,J. and Swinnerton-Dyer, Sir Peter: The Hasse problem for rational surfaces, in: Collection of articles dedicated to Helmut Hasse on his seventy-fifth birthday~III, {\em J.\ Reine Angew.\ Math.} {\bf 274/275}\br(1975)\brr164--174
\bibitem[BMS]{BMS}
Browning, T.\,D., Matthiesen, L., and Skorobogatov, A.\,N.: Rational points on pencils of conics and quadrics with many degenerate fibers, {\em Ann.\ of Math.} {\bf 180}\br(2014)\brr381--402
\bibitem[Ch]{Ch}
Chrystal, G.: Algebra, An elementary text-book for the higher classes of secondary schools and for colleges, Part~II, Reprint of the 6th edition, {\em Chelsea Publishing Co.}, New York~1959
\bibitem[De]{De}
Demazure, M.: Surfaces de Del Pezzo II--\smash{\'Eclater}
$n$
points dans
$\Pb^2$,
{\em S\'eminaire sur les singularit\'es des surfaces (Polytechnique)} (1976--1977), Expos\'e n$^\circ$4, 1--13
\bibitem[Do]{Do}
Dolgachev, I.\,V.: Classical Algebraic Geometry: a modern view, {\em Cambridge University press,} Cambridge~2012
\bibitem[Ha]{Ha}
Hartshorne, R.: Algebraic Geometry, Graduate Texts in Mathematics 52, {\em Springer,} New York 1977
\bibitem[Is]{Is}
Iskovskih, V.\,A.: Minimal models of rational surfaces over arbitrary fields (Russian), {\em Izv.\ Akad.\ Nauk SSSR\/} {\bf 43}\br(1979)\brr19--43
\bibitem[Ja]{Ja}
Jahnel, J.: Brauer groups, Tamagawa measures, and rational points on algebraic varieties, Mathematical Surveys and Monographs~198, {\em AMS,} Providence~2014
\bibitem[Man]{Man}
Manin, Yu.\,I.: Cubic forms, algebra, geometry, arithmetic,
{\em North-Holland Publishing Co.} and {\em American Elsevier Publishing Co.,}
Amsterdam, London, and New York~1974
\bibitem[Maz]{Maz}
Mazur, B.: The topology of rational points, {\em Experiment.\ Math.} {\bf 1}\br(1992)\brr35--45
\bibitem[Mi]{Mi}
Milne, J.\,S.: \smash{\'Etale} Cohomology, {\em Princeton University Press,} Princeton~1980
\bibitem[Ne]{Ne}
Neukirch, J.: Algebraic number theory, Grundlehren der Mathematischen Wis\-sen\-schaf\-ten 322, {\em Springer,} Berlin~1999
\bibitem[Pi]{Pi}
Pierce, R.\,S.: Associative algebras, Graduate Texts in Mathematics~88, {\em Springer,} New York-Berlin~1982
\bibitem[SSk]{SSk}
Salberger, P.\ and Skorobogatov, A.\,N.: Weak approximation for surfaces defined by two quadratic forms, {\em Duke Math.\ J.} {\bf 63}\br(1991)\brr517--536
\bibitem[Sch]{Sch}
Schl\"af\/li,~L.: An attempt to determine the twenty-seven lines upon a surface of the third order, and to divide such surfaces into species in reference to the reality of the lines upon the surface, {\em Quart.\ J.\ Math.} {\bf 2}\br(1858)\brr110--120
\bibitem[Se]{Se}
Serre, J.-P.: Cours d'arithm\'etique, {\em Presses Universitaires de France,} Paris~1977
\bibitem[Silh]{Silh}
Silhol, R.: Real algebraic surfaces, Lecture Notes in Mathematics~1392, {\em Springer,} Berlin~1989
\bibitem[Silv]{Silv}
Silverman, J.\,H.: The arithmetic of elliptic curves, Second edition, Graduate Texts in Mathematics~106, {\em Springer,} Dordrecht~2009
\bibitem[Sk]{Sk}
Skorobogatov, A.\,N.: Torsors and rational points, Cambridge Tracts in Mathematics~144, {\em Cambridge University Press,} Cambridge~2001
\bibitem[SD1]{SD1}
Swinnerton-Dyer, Sir~Peter: Two special cubic surfaces, {\em Mathematika\/} {\bf 9}\br(1962)\brr54--56
\bibitem[SD2]{SD2}
Swinnerton-Dyer, Sir~Peter: The Brauer group of cubic surfaces, {\em Math.\ Proc.\ Cambridge Philos.\ Soc.} {\bf 113}\br(1993)\brr449--460
\bibitem[Ta]{Ta}
Tate,~J.: Global class field theory, in: Algebraic number theory, Edited by J.\,W.\,S.~Cassels and A.~Fr\"ohlich, {\em Academic Press\/} and {\em Thompson Book Co.,} London and Washington~1967, 162--203
\bibitem[VA]{VA}
V\'arilly-Alvarado, A.: Arithmetic of del Pezzo surfaces, Notes of lectures given at the Lorentz Center, Leiden, in October~2010, available at: {\tt http://math.rice.edu/\~{}av15/Files/\discretionary{}{}{}LeidenLectures.pdf}
\bibitem[VAV]{VAV}
V\'arilly-Alvarado, A.\ and Viray, B.: Arithmetic of del Pezzo surfaces of degree 4 and vertical Brauer groups, {\em Adv.\ Math.} {\bf 255}\br(2014)\brr153--181
\bibitem[Wi]{Wi}
Wittenberg, O.: Intersections de deux quadriques et pinceaux de courbes de
genre~$1$,
Lecture Notes in Mathematics~1901, {\em Springer,} Berlin~2007
\end{thebibliography}
\end{document}